\documentclass[11pt]{article}

\usepackage{amsmath}
\usepackage{amssymb}

\title{THE PROPAGATION OF MOLECULAR CHAOS BY MARKOV TRANSITIONS}

\author{By Alexander Gottlieb}

\date{ }

\begin{document}

\def \e {\varepsilon}
\def \pa {\partial} 

\def \bx{{\bf x}}
\def \by{{\bf y}}
\def \bv{{\bf v}}
\def \bw{{\bf w}}
\def \bl{{\bf l}}
\def \bs{{\bf s}}
\def \bt{{\bf t}}
\def \bxp{{\bf x'}}
\def \byp{{\bf y'}}
\def \bvp{{\bf v'}}
\def \bwp{{\bf w'}}
\def \bvs{{\bf v^*}}
\def \bws{{\bf w^*}}

\def \ddt{ \frac{d}{dt}}
\def \ppt{ \frac{\partial}{\partial t}}
\def \oon {\frac{1}{n}}
\def \EE {\mathbb{E}}
\def \PP {\mathbb{P}}
\def \RR {\mathbb{R}}
\def \NN {\mathbb{N}}
\def \RRd {\mathbb{R}^d}

\def \B {\mathcal{B}}
\def \F {\mathcal{F}}
\def \G {\mathcal{G}}
\def \I {\mathcal{I}}
\def \J {\mathcal{J}}
\def \P {\mathcal{P}}
\def \T {\mathcal{T}}
\def \Kt {\widetilde{K}}
\def \phat {\widehat{\phi}}

\def \s {\sigma}

\newtheorem{proposition}{Proposition}[section]
\newtheorem{theorem}{Theorem}[section]
\newtheorem{lemma}{Lemma}[section]
\newtheorem{definition}{Definition}[section]
\newtheorem{corollary}{Corollary}[section]

\maketitle

\begin{abstract}

We establish a necessary and sufficient condition for the propagation
of chaos by a family of many-particle
Markov processes, if the particles live in a Polish space $S$:  A
sequence of $n$-particle Markov transition functions $\{K_n\}$ propagates
molecular chaos if and only if the sequence $\{K_n(\bs_n,\cdot)\}$ is
chaotic whenever $\bs^n = (s^n_1,s^n_2,\ldots,s^n_n) \in S^n$ is such that 
$\oon \sum_i \delta(s^n_i)$
converges to a law on $S$ as $n \longrightarrow \infty$.
\footnote{This article is a short version of the author's doctoral
dissertation \cite{Gottlieb}.  The work was supported in part by the Office
of Science, Applied Mathematical Sciences Subprogram, of the
U.S. Department of Energy, under Contract Number DE-AC03-76SF0098.}
\footnote{ 
{\it AMS} 1981 {\it subject classifications}.  Primary 82C05.
Secondary 60J35, 60K35.

 }
\end{abstract}

\section{Introduction}

  The ``propagation of chaos'' means the persistence 
in time of molecular chaos (i.e., the stochastic independence of two random particles
in a many-particle system) in the limit of infinite particle number.
Propogation of chaos is an important concept of kinetic theory that
relates the equations of Boltzmann and Vlasov to the dynamics
of many-particle systems.

The concept of propagation of chaos originated with Kac's Markovian models
of gas dynamics \cite{Kac55,Kac}.   Kac invented a class of interacting particle
systems wherein particles collide at random with
each other while the {\it density} of particles evolves
deterministically in the limit of infinite particle number.  
A nonlinear evolution equation analogous to Boltzmann's
equation governs the particle density.
The processes of Kac were further investigated with regard to
their fluctuations about the deterministic infinite particle limit 
in \cite{McK75, Uch83, Uch88}. 
McKean introduced propagation
of chaos for interacting diffusions and analyzed what are now called
McKean-Vlasov equations \cite{McK66, McK68}.  Independently,
Braun and Hepp \cite{BH} analyzed the propagation of chaos for Vlasov
equations and proved a central limit theorem for the fluctuations.
Analysis of the fluctuations and large deviations for McKean-Vlasov
processes was carried out in \cite{Tanaka, Sz84, DawsonGartner}.
Other instances of the propagation of chaos have been studied in
\cite{Oelsch, RezTar, Graham}.  
Finally, we refer the reader to the long, informative articles by
 Sznitman~\cite{Sznitman} and by
M\'el\'eard~\cite{Meleard} in Springer-Verlag's {\it Lecture Notes in
Mathematics}.  
 
This article is organized as follows:  

In Section \ref{Molecular Chaos}, molecular chaos is defined and a class
of examples is provided.  The
weak convergence approach to the study of molecular chaos due to Sznitman and Tanaka
is introduced in Theorem \ref{SznTan}.   In Section \ref{The
Propagation of Chaos}, propagation of chaos is defined, and the main
results of this article are stated in Theorem \ref{main} and its
corollary.  Theorem \ref{main} is proved in
Section \ref{Proof Section} using Prohorov's theorem and Theorem
\ref{SznTan} of Section \ref{Molecular Chaos}.

\section{Molecular Chaos}
\label{Molecular Chaos}

In what follows, if $(S,d_S)$ is a separable metric space, its Borel
algebra will be denoted $\B(S)$, the space of bounded and continuous
functions on $S$ will be denoted $C(S)$, and the space of probability
laws on $S$ with the weak topology will be denoted $\P(S)$.

Let $S^n$ denote the n-fold product of $S$ with itself; 
\[
      S^n := \left\{ (s_1, s_2, \ldots, s_n): s_i \in S \quad \mathrm{ for} 
       \quad i=1,2,\ldots,n \right\}  .
\] 
$S^n$ is itself metrizable in a variety of equivalent ways that all 
generate the same topology and the same Borel algebra $\B({S^n})$. 
If $\rho_n$ is a law on $S^n$ and $k \le n$, let $\rho_n^{(k)}$ denote
the marginal law on the first $k$ coordinates, that is, the law on
$S^k$ induced by the map $(s_1,s_2,\ldots,s_n) \mapsto
(s_1,\ldots,s_k)$.  If $\rho \in \P(S)$, let $\rho^{\otimes n}$ denote
the product law on $S^n$ defined by $\int_{S^n} f_1(s_1)\cdots
f_n(s_n) \rho^{\otimes n}(ds_1 \cdots ds_n) = \prod_{i=1}^n \int_S
f_i(x)\rho (dx)$.

Let $\Pi_n$ denote the set of permutations of $\{1,2,\ldots,n\}$.   
The permutations $\Pi_n$ act on $S^n$ by permuting coordinates:
the map $\pi\cdot:S^n \longrightarrow S^n$ is
\[
    \pi \cdot (s_1, s_2, \ldots, s_n)  := (s_{\pi(1)},s_{\pi(2)}, 
    \ldots,s_{\pi(n)})  .
\]
If $E$ is any subset of $S^n$, define 
\[
     \pi\cdot E = \{ \pi\cdot\bs: \bs \in E \}  .
\]
A law $\rho$ on $S^n$ is {\it symmetric} if $\rho(\pi\cdot B) = \rho(B)$ 
for all $\pi \in \Pi_n$ and all $B \in \B({S^n})$.  Products
$\rho^{\otimes n}$ are symmetric, for example.  The {\it 
symmetrization} ${\widetilde \rho}$ of a law $\rho \in \P(S^n)$ is the 
symmetric law such that 
\[
    {\widetilde \rho}(B) := \frac{1}{n!}\sum_{\pi \in \Pi_n} \rho(\pi\cdot 
    B),
\]
for all $B \in \B({S^n})$. 

For each $n$, define the map $\e_n:S^n \longrightarrow \P(S)$ by  
\begin{equation}
\label{e's}
        \e_n((s_1,s_2,\ldots,s_n)) := \oon\sum_{i=1}^n \delta(s_i),
\end{equation} 
where $\delta(x)$ denotes a point-mass at $x$.
These maps are measurable for 
each $n$, and $\e_n(\pi \cdot\bs)  = \e_n(\bs)$ for all $\bs \in S^n, \pi \in \Pi_n$.  
               
To derive the Boltzmann equation for a dilute gas, Boltzmann assumed
that a condition of ``molecular disorder'' obtains, i.e., the velocities
of two random molecules in a gas are stochastically independent \cite{Boltzmann}.
This is not an entirely realistic assumption: surely the collisions
between the molecules of an $n$-particle gas must introduce some
stochastic dependence, even if the molecules were independent to begin
with.  Nonetheless, if the number $n$ of molecules is very large, one
might well expect that the velocities of two random molecules would be
{\it nearly} independent; the presence of so many other molecules
should wash out most of the dependence between two randomly selected
molecules.  These considerations motivated Kac to define molecular
disorder as a kind of asymptotic independence of pairs of particles in
the infinite-particle limit.  He called the definiendum the
``Boltzmann property'' instead of ``molecular chaos'' or simply
``chaos,'' as it is now termed.  Kac's definition is

\begin{definition}[Kac, 1954]
\label{chaos}
Let $(S,d_S)$ be a separable metric space.  Let $\rho$ be a law on $S$, and for $n = 1,2 ,\ldots$,
let $\rho_n$ be a symmetric law on $S^n$.  

The sequence
$\{ \rho_n \}$ is $\rho$-{\bf chaotic} if, for each natural number $k$
and each choice $   
     \phi_1(s), \phi_2(s), \ldots, \phi_k(s)
$
 of $k$ bounded and continuous functions on $S$,
\[
\lim_{n \rightarrow \infty} \int_{S^n} \phi_1(s_1) \phi_2(s_2) \cdots
\phi_k(s_k) \rho_n(ds_1 ds_2\ldots ds_n)    =
\prod_{i=1}^k \int_S \phi_i(s) \rho(ds).
\]
\end{definition}

It is now known \cite{Szn, Tanaka, Gottlieb} that a sequence of symmetric laws
$\{ \rho_n \}$ (with $\rho_n \in \P(S^n)$) is $\rho$-chaotic if and
only if the marginals $\rho_n^{(k)}$ converge weakly to $\rho^{\otimes
k}$ as $n \longrightarrow \infty$ for each fixed $k$.  (We suggest that
this last condition be taken for the {\it definition} of chaos in
the future, in view of its simplicity and equivalence to Kac's
definition.)  The following theorem gives conditions that are
equivalent to molecular chaos.

\begin{theorem}[Sznitman, Tanaka]
\label{SznTan}
Let $S$ be a separable metric space, and for each $n$ let $\rho_n$ be 
a symmetric law on $S^n$.  

The following are equivalent:
\begin{trivlist}
\item{(i)} $\{ \rho_n \}$ is $\rho$-chaotic;
\item{(ii)}  For all $\phi_1,\phi_2 \in C_b(S)$,
\[
  \lim_{n \rightarrow \infty} \int_{S^n} \phi_1(s_1) \phi_2(s_2) 
  \rho_n(d\bs)  = \int_S \phi_1(s) \rho(ds) \int_S \phi_2(s) \rho(ds) ;
\]
\item{(iii)} The marginals $\rho_n^{(2)}$ converge weakly to $\rho\otimes\rho$ 
  as $n$ tends to infinity.
\item{(iv)}  For all $k$, 
  the marginals $\rho_n^{(k)}$ converge weakly to $\rho^{\otimes k}$ 
  as $n$ tends to infinity.
\item{(v)}  The laws $\rho_n\circ 
  \e_n^{-1}$ converge to $\delta(\rho)$ in $\P(\P(S))$ as $n$ tends to
  infinity;

\end{trivlist} 

\end{theorem}

An important class of chaotic sequences arises in connection with
 the principle of equivalence of microcanonical and canonical
averaging of statistical mechanics.  The principle follows from
 the molecular chaos of sequences of microcanonical distributions.
For example, let $S$ be a finite set and let $H:S
\longrightarrow \RR$ be any function that we call energy:
$H(s)$ is the energy of a particle if it is the state $s \in S$.  Fix 
$E \in \RR$, to be interpreted as average energy per particle, and
$\delta > 0$, and define $\mu_{n,E,\delta} \in \P(S^n)$ to be the law
that gives equal probability to each $n$-particle state
$(s_1,s_2,\ldots,s_n)$ that satisfies $\oon \sum_{i=1}^n H(s_i) \in (E
- \frac{1}{2} \delta, E + \frac{1}{2}\delta)$ and probability $0$ to all other points of
$S^n$.  Thus, the ``microcanonical distributions'' $\mu_{n,E,\delta}$
render equiprobable all $n$-particle
states for which the energy-per-particle is approximately $E$.  The
sequence of microcanonical distributions $\{ \mu_{n,E,\delta} \}_{n
\in \NN}$ is $\gamma$-chaotic, where $\gamma \in \P(S)$ is the
Gibbsian distribution defined by 
\[
    \gamma(t) = e^{-\beta H(t)} / \sum_{s \in S} e^{-\beta H(s)}
\]
for all $t \in S$, where the parameter $\beta$ depends on $E$ and
$\delta$.  This follows from Sanov's large deviations theorem.

In closing, we note the following pleasant corollary of Theorem
\ref{SznTan} (see \cite{Gottlieb} for a proof):
\begin{corollary}
\label{SpecEntThm}
Let $S = \{s_1,s_2,\ldots,s_k\}$ be a finite set, 
and for each $n$ let $\rho_n \in \P(S^n)$ be a symmetric law.

If the sequence $\{\rho_n\}$ is $ p$-chaotic, 
then 
\[
  \lim_{n \rightarrow \infty} \left( \oon \sum_{{\bf x}\in S^n} 
      \rho_n({\bf x}) \log \rho_n({\bf x}) \right) =
      \sum_{i=1}^k p_i \log p_i.
\]
\end{corollary}

\section{The Propagation of Chaos}
\label{The Propagation of Chaos}

Let $(S, d_S)$ and $(T,d_T)$ be two separable metric spaces.
A Markov transition function $K$ is a 
function on $S \times \B(T)$ that satisfies the following two conditions:

\noindent (1) \qquad $K(s,\cdot)$ is a probability 
measure on $(T,\B(T))$ for each $s \in S$, and

\noindent (2) \qquad $K(\cdot, E)$ is a measurable function on 
$(S,\B(S))$ for each $E \in \B(T)$.

A Markov transition function on $S\times \B(T)$ is also
called a Markov transition from $S$ to $T$.

A Markov {\it process} on a state space $(S,d_S)$ determines a 
family, indexed by time, of Markov transitions from $S$ to
itself: $\{ K(s,E,t) \}_{t \ge 0}$.  
The transitions 
satisfy the Chapman-Kolmogorov equations 
\[
      K(s,E,t + t') = \int_S K(s,dx,t) K(x,E,t'); \quad t,t' \ge 0,s \in S, E \in \B(S)
\]
in addition to (1) and (2) above.

For each $n$, let $K_n$ be a Markov 
transition from $S^n$ to $T^n$.  
We assume that the Markov transition function $K_n$ 
is such that 
if $\pi$ is a permutation in $\Pi_n$  and $A$ is a Borel 
subset of $T^n$, 
\begin{equation}
    K_n(\pi\cdot \bs,\pi\cdot A) =  K_n( \bs, A)   .   
\label{permute}
\end{equation}

\begin{definition}
\label{PoCDef}

A sequence $\{K_n\}$ satisfying (\ref{permute}) {\bf propagates chaos} if,
 whenever $\{ \rho_n \}$ is a $\rho$-chaotic sequence of 
symmetric laws on $S^n$, the sequence
\[
  \left\{ \int_{S_n} K_n({\bf s},\cdot) \rho_n(d{\bf s}) \right\}_{n=1}^{\infty}
\]
is $\tau$-chaotic for some $\tau \in \P(T)$.
\end{definition}

 When we say that a family of $n$-particle Markov {\it processes} on a state space $S$
{\bf  propagates chaos} we mean that, for each fixed time $t > 0$, the family of associated
 $n$-particle transition functions $\{K_n(\bs,E,t) \}$ 
 propagates chaos.  Propagation of chaos will be seen to imply the
existence of a semigroup of (typically nonlinear) operators $ F_t$ 
on $\P(S)$ such that 
$\{ \int K_n(\bs,\cdot)\rho_n(d\bs) \}$ is $F_t\rho$-chaotic if $\{ \rho_n \}$
is $\rho$-chaotic. 
For families of interacting particle systems suited 
to the study of gases or plasmas, the semigroup 
$\{ F_t \}_{t \ge 0} $ is
 the semigroup of solution operators for the Boltzmann or the Vlasov 
equation. 

In case $S$ is Polish, and if Definition \ref{PoCDef} for the
propagation of chaos is accepted (see Remark 2 at the end of this
section), we have a necessary and sufficient
condition for chaos to be propagated, which is expressed in terms of the maps
$\e_n$ defined in (\ref{e's}).

\begin{theorem}
\label{main}
Suppose $(S, d_S)$ is a complete, separable metric space.  Let
$\{K_n\}$ be a sequence of Markov transitions that satisfy condition
(\ref{permute}), and denote the symmetrization of $K_n(\bs_n,\cdot)$
by $\Kt_n(\bs_n,\cdot)$. 

$\{ K_n \}$ 
propagates chaos if and only if there exists a continuous function
\[ F:\P(S) \longrightarrow \P(T)  \]
such that 
$     \left\{ \Kt_n(\bs_n,\cdot) \right\} $
is $F(p)$-chaotic whenever the points $\bs_n \in S^n$ are such that 
$ \e_n(\bs_n) \longrightarrow p $ in $\P(S)$.   

\end{theorem}

Propagation of chaos is easily characterized in case the
Markov transitions are deterministic, that is, when the $n$-particle dynamics are
given by measurable maps from $S^n$ to $T^n$ that
commute with permutations of coordinates.
Let $k_n:S^n \longrightarrow T^n$ be a measurable map
that commutes with permutations of $n$-coordinates, i.e., such that
\[
       k_n(s_{\pi(1)},s_{\pi(2)},\ldots,s_{\pi(n)}) = \pi \cdot 
               k_n(s_1,s_2,\ldots,s_n)
\]
for each point ${\bf s} \in S^n$ and each
permutation $\pi$ of the symbols ${1,2,\ldots,n}$.
Given $k_n$, define the Markov transition $K_n$ from $S^n$ to $T^n$ by
\[
       K_n({\bf s},E) = {\bf 1}_E(k_n({\bf s}))
\]
when ${\bf s} \in S^n$ and $E \in \mathcal{B}_{T^n}$.
Say that $\{ k_n \}_{n=1}^{\infty}$ {\it propagates chaos} if the 
sequence of deterministic 
transition functions $\{ K_n \}$ propagates chaos.

\begin{corollary}[Deterministic Case]
\label{deterministic}
Let $(S,d_S)$ be a complete and separable metric space, and for each $n$ let  
$k_n$ be a measurable map from $S^n$ to $T^n$ that commutes with 
permutations.

$\left\{ k_n \right\} $ propagates chaos if and only if there exists a continuous function    
\[
      F:\P(S) \longrightarrow \P(T) 
\]
such that $  \e_n(k_n({\bf s}_n)) \longrightarrow F(p) $
in $\P(T)$ whenever $ \e_n({\bf s}_n) \longrightarrow p $
in $\P(S)$.
\end{corollary}

\noindent {\bf Remark 1}: \qquad
Most Markov processes of interest are characterized by their laws on nice 
 {\it path spaces}, such as $C([0,\infty), S)$ the space of continuous
paths in $S$, or the space $D([0,\infty), S)$ of right continuous paths in $S$
 having left limits.  For such processes, 
 the function that maps a state $s \in S$ to the law of the process 
 started at $s$ defines a Markov transition from $S$ to the entire 
 path space.  
Now, if a sequence of transitions $K_n(s,\cdot)$ from $S^n$ 
 to the path spaces $C([0,\infty),S^n)$ or $D([0,\infty),S^n)$ propagates 
 chaos, then, {\it a fortiori}, it propagates chaotic sequences of
initial laws to chaotic sequences of laws on $S^n$ at any (fixed) later
time.  We have defined the propagation of chaos
for sequences of Markov transitions from $S^n$ to a (possibly) different space $T^n$, instead of
simply from $S^n$ to itself, with the case where $T$ is path space
especially in mind.
This way, our study will pertain even
to those families of processes that propagate the chaos of initial
laws to the chaos of laws on the whole path space.

\noindent {\bf Remark 2}: \qquad
We are adopting here a strong definition of the propagation of chaos.
Other authors \cite[p. 42]{Meleard}\cite[p. 98]{Pulvirenti}
define propagation of chaos in a {\it weaker}
sense (before going to prove that chaos propagates in a variety of
important cases).  According to these authors, 
a family of Markovian $n$-particle processes propagates chaos
if  $\{\int K_n(\bs,\cdot) \rho^{\otimes n}(d\bs)\}$ is chaotic for all
$\rho \in \P(S)$  and $t > 0$, where $ \rho^{\otimes n}$ is 
product measure on $S^n$.   In other words, 
only {\it purely} chaotic sequences of initial measures
are required to ``propagate'' to chaotic sequences.   This condition 
is strictly weaker than the one we adopt for our definition.  
For example, take $S = \{0,1\}$ and let $\delta(x)$ or $\delta_x$
denote a point mass at $x$.  Then, if 
\[
   K_n(\bs,\cdot,t) = \left\{ \begin{array}{cc}
                 \delta_{(1,1,\ldots,1)} & \mathrm{if}\  \bs \ne (0,0,\ldots,0) \\
                 \delta_{(0,0,\ldots,0)} & \mathrm{if}\  \bs = (0,0,\ldots,0)   \\
                  \end{array} \right.
\]
for all $t>0$, the sequence $\{K_n\}$ propagates chaos in the weak sense, but
not in the strong sense of our definition.   Under these $K_n$'s,
 the $\delta(0)$-chaotic sequence $\{ \delta_{(0,0,\ldots,0)} \}$ is
propagated to itself, while other $\delta(0)$-chaotic sequences
are propagated to $\delta(1)$-chaotic sequences, and yet other
$\delta(0)$-chaotic sequences are not propagated to chaotic sequences
at all.

\section{Proof of Theorem \ref{main} }
\label{Proof Section}

Let $(X,d_X)$ be a metric  space, and $D_1 \subset D_2 \subset \cdots $ an 
increasing chain of Borel subsets of $X$ whose union is dense in $X$.  For each 
natural number $n$, let $f_n$ be a measurable real-valued function on $D_n$.

Consider the following four conditions on the sequence $\{ f_n \}$.  
They are listed in order of decreasing 
strength.

\begin{trivlist} 

\item {[A]} \qquad Whenever $\{ \mu_n \}$ is a weakly convergent sequence of
probability measures on $X$ with $\mu_n$ supported on $D_n$, then 
the sequence 
\[ 
\left\{ \int_X f_n(x) \mu_n(dx) \right\}_{n=1}^{\infty} 
\]
of real numbers converges as well.

\item {[B]} \qquad Whenever $\{ \mu_n \}$ is a sequence of 
probability measures on $X$ that converges weakly to $\delta(x)$ for some 
$x \in X$, and $\mu_n$ is supported on $D_n$, then the 
sequence $\left\{ \int_X f_n(x) \mu_n(dx) \right\}$ also converges.

\item {[C]} \qquad Whenever $\{ d_n \}$ is a convergent 
sequence of points in $X$, with $d_n \in D_n$, then $\{ f_n(d_n) \}$ also 
converges.

\item {[D]} \qquad There exists a continuous function $f$ on $X$ 
towards which the functions $f_n$ converge uniformly on compact sets: 
for any compact $K \subset X$, and for any $\epsilon > 0$, 
there exists a natural number $N$ such that, whenever $n \ge N$ and 
$d \in D_n \cap K$, then $|f(d) - f_n(d)| < \epsilon  . $

\end{trivlist}

\begin{lemma}
\label{lemma}
If $X$ is a Polish space, i.e., if $X$ is homeomorphic to a complete and
separable metric space, and $\sup_{d \in D_n}\{ |f_n(d)| \} \le B$ for all $n$, 
then conditions [A], [B], [C], and [D] are all equivalent.
\end{lemma}
\noindent {\bf Sketch of Proof}:

[A] $\implies$ [B] $\implies$ [C] $\implies$
[D] even when $X$ is not Polish.  

 To get [D] $\implies$ [A], we are
assuming that $X$ is Polish, for then there must exist, by Prohorov's
theorem, compact sets which support all of
the measures of the convergent sequence $\{ \mu_n \}$ to within any
$\epsilon > 0$.  (Prohorov's theorem states that
if $(X,\T)$ is Polish, then $\Sigma \subset \P(X)$ is tight if and
only if its closure is compact in $\P(X)$.)
$\qquad \blacksquare$

Let $S$ and $T$ be Polish spaces.  For each natural 
number $n$, 
let $K_n$ be a Markov transition from $S^n$ to $T^n$ that satisfies 
the permutation condition (\ref {permute}). 
Henceforth, the notation $\e_n$  
is used both for the map from $S^n$ to $n$-point empirical measures on 
$S$ and for the same kind of map on $T^n$.   $\e_n(S^n)$ denotes the
subset of $\P(S)$ that consists of discrete laws of the form
(\ref{e's}).

Now, Markov transitions $K_n$ from $S^n$ to $T^n$ induce 
Markov transition functions $H_n$ from $\e_n(S^n)$ to $\e_n(T^n)$ defined by
\begin{equation}
  H_n(\zeta,G) := \int_{{\bf s} \in S^n} J_n(\zeta, d{\bf s}) 
        K_n\left( {\bf s}, \e_n^{-1} ( G ) \right)  
\label{InducedTransition}
\end{equation} 
for $\zeta \in \e_n(S^n)$ and $G$ a measurable subset of  $\e_n(T^n)$,
where $J_n$ is the Markov transition from $\e_n(S^n)$ to $S^n$ such
that $J_n(\zeta,\cdot)$ that allots equal probability to 
each of the points in $\e^{-1}(\{ \zeta \})$.  Note that $J_n(\zeta,\cdot) \circ
\e_n^{-1} = \delta(\zeta)$.
      
Theorem \ref{SznTan} shows that propagation of chaos by a 
sequence $\{ K_n \}$ is equivalent to the following condition on the
sequence of induced transitions $\{H_n\}$:
\begin{lemma}
\label{PoCDef2}
The sequence of Markov transitions $\{ K_n \}$ propagates chaos if 
and only if, 
whenever $\{ \mu_n \} $  converges to 
$\delta(p) \in \P(\P(S))$ for some $p \in \P(S)$, 
with $\mu_n$ supported on $\e_n(S_n)$ for each $n$, the sequence
\begin{equation}
\label{H sequence} 
\left\{ \int_{\e_n(S^n)} H_n(\zeta,\cdot) \mu_n(d\zeta) \right\}_{n=1}^\infty
\end{equation}
converges in $\P(\P(T))$ to $\delta(q)$, for some $q \in \P(T)$.
\end{lemma}
\noindent{\bf Proof}:

By condition (v) of Theorem \ref{SznTan}, $\{ \mu_n \in \P(\e_n(S^n)) \}$ converges to  
a point mass $\delta(p) \in \P(\P(S))$ if and only if 
$
   \left\{ \int_{\e_n(S^n)} J_n(\zeta,\cdot) \mu_n(d\zeta) 
   \right\}
$
is chaotic.  Therefore, the sequence of transitions $\{ 
K_n \}$ propagates chaos if and only if 
\begin{equation}
\label{K sequence}
   \left\{  \int_{S^n} 
   K_n(\bs,\cdot) \int_{\e_n(S^n)} J_n(\zeta,d\bs) \mu_n(d\zeta) 
   \right\}_{n=1}^{\infty}
\end{equation} 
is chaotic whenever $ \mu_n \longrightarrow \delta(p)$.  
By part (v) of Theorem \ref{SznTan} again,
 the sequence (\ref{K sequence}) is chaotic if and only if 
\[
 \left( \int_{S^n}  K_n(\bs,\cdot) \int_{\e_n(S^n)} J_n(\zeta,d\bs)
\mu_n(d\zeta) \right) \circ \e_n^{-1} \longrightarrow \delta(q)
\]
 for some $q \in \P(T)$.
But, by definition (\ref{InducedTransition}) of the transitions $H_n$,
\[
   \left( \int_{S^n} 
     K_n(\bs,\cdot) \int_{\e_n(S^n)} J_n(\zeta,d\bs) \mu_n(d\zeta)   
     \right) \circ \e_n^{-1}  =
    \int_{\e_n(S^n)}  H_n(\zeta,\cdot)  
   \mu_n(d\zeta),  
\]
so $\{K_n\}$ propagates chaos if and only if the sequence
(\ref{H sequence}) converges to $\delta(q)$ for some $q \in \P(T)$ whenever
$\{\mu_n\}$ converges to $\delta(p) \in \P(\P(S))$.
\qquad $\blacksquare$

Having Lemmas \ref{lemma} and \ref{PoCDef2} in hand, we proceed to the

\medskip
\noindent {\bf Proof of Theorem \ref{main}}:

  For each bounded and continuous function $\phi \in C_b(\P(T))$ 
define the functions $ \phat_n: \e_n(S^n) \longrightarrow {\mathbb R}$ by
\begin{equation}
\label{phat}
     \phat_n(\zeta) := \int_{\P(T)} \phi(\eta) H_n(\zeta,d\eta) , 
\end{equation}
 where $H_n$ is as defined in (\ref{InducedTransition}).  Note that
the functions $\phat_n $ are bounded uniformly in $n$ since $\phi$ is bounded.

Assume first that $\{K_n\}$ propagates chaos.  We will show that there
exists a continuous function $F:\P(S) \longrightarrow \P(T)$ such that
$  \left\{ \Kt_n(\bs_n,\cdot) \right\} $
is $F(p)$-chaotic whenever the points $\bs_n \in S^n$ are such that 
$ \e_n(\bs_n) \longrightarrow p $ in $\P(S)$.

Since $\{K_n\}$ propagates chaos, Lemma \ref{PoCDef2} implies that
whenever $\{ \mu_n \in \P(\e_n(S^n)) \}$ converges in $\P(\P(S))$ to 
$\delta(p)$, then
\begin{equation}
\label{thingy} 
   \int_{\e_n(S^n)} \phat_n(\zeta) \mu_n(d\zeta) \longrightarrow \phi(q)
\end{equation}
for some $q \in \P(T)$.  In fact, Lemma \ref{PoCDef2} implies 
that $q$ does not depend on our choice of $\phi$: the same $q$ works
for all $\phi$ in (\ref{thingy}).  

Condition (\ref{thingy}) resembles condition [B] above, so Lemma
\ref{lemma} shows that there exists a
continuous function $G_{\phi}(p)$, depending on $\phi$, such that if 
 $\{ \bs_n \in S^n \}$ is a sequence satisfying
$ \e_n(\bs_n) \longrightarrow p $ in $\P(S)$, then 
\[
     \phat_n(\e_n(\bs_n))  \longrightarrow G_{\phi}(p).
\]
By (\ref{thingy}), $G_{\phi}(p) = \phi(q)$ for some $q \in \P(T)$ {\it
that does not depend on} $\phi$.  The only way that all the
$G_{\phi}$'s can have this form and yet all be continuous is for the
dependence of $q$ on $p$ to be continuous: there
must be a continuous $F$ from $\P(S)$ to $\P(T)$ such that
$G_{\phi}(p) = \phi(F(p))$ for all $\phi \in C_b(\P(T))$.

Thus, there exists a continuous function $F$ from $\P(S)$ to $\P(T)$
such that 
\[
    \left[ \e_n(\bs_n) \longrightarrow p \right] \implies 
    \left[ \phat_n(\e_n(\bs_n))  \longrightarrow \phi(F(p)) \right]
\]
for all $\phi \in C_b(\P(T))$.  This fact, and the definitions
(\ref{InducedTransition}) and (\ref{phat}) of $H_n$ and $\phat$, imply
that 
\[
    \left[ \e_n(\bs_n) \longrightarrow p \right] \implies 
    \left[ \Kt_n(\bs_n, \cdot) \circ \e_n^{-1} \longrightarrow \delta(F(p)) \right].
\]
Finally, by Theorem \ref{SznTan}, we have that 
\[
    \left[ \e_n(\bs_n) \longrightarrow p \right] \implies 
    \left\{ \Kt_n(\bs_n, \cdot)\right\}_{n=1}^{\infty} \  \mathrm{is} \  
     F(p)\mathrm{-chaotic} . 
\]

This demonstrates the necessity of the condition of Theorem
\ref{main}.  
Next we demonstrate its sufficiency, i.e., that $\{ K_n \} $
propagates chaos {\it if} there exists a
continuous function $F:\P(S) \longrightarrow \P(T)$ such that 
$\{ \Kt_n(\bs_n,\cdot) \}$ is $F(p)$-chaotic whenever 
$\e_n(\bs_n)$ converges weakly to $p \in \P(S)$.

Suppose $P_n \in \P(S^n)$ is $p$-chaotic.  Let $\mu_n = P_n \circ \e_n^{-1}$.  
Then $\mu_n \longrightarrow \delta(p)$ in $\P(\P(S))$ by Theorem 
\ref{SznTan}.   Our goal is to prove that 
\[
    \int_{\P(S)} H_n(\zeta, d\eta) \mu_n(d\zeta) \longrightarrow 
    \delta(F(p)),    
\]
where $H_n$ is as defined in (\ref{InducedTransition}).
This is enough, by Lemma \ref{PoCDef2}, to demonstrate that 
chaos propagates.

By hypothesis, if $\{\bs_n \in S^n\}$ is such that $\e_n(\bs_n)$ 
converges to $p$ then $\{ \Kt_n (\bs_n,\cdot) \}$ is $F(p)$-chaotic.  
By Theorem \ref{SznTan} and the fact that
\[ 
    \Kt_n(\bs_n, \cdot)\circ \e_n^{-1} = H_n(\e_n(\bs_n),\cdot) ,
\]
the hypothesis is equivalent to the statement that, if
$ p_n \in \e_n(S^n)$ for each $n$, then
\begin{equation}
\label{hypo2}
     \left[ p_n \longrightarrow p \right] 
     \implies \left[H_n(p_n,\cdot) \longrightarrow \delta(F(p)) 
     \right].  
\end{equation}

The hypothesis (\ref{hypo2}) and equation (\ref{phat}) imply that for
each $\phi \in C_b(\P(T)) $,
\[
   \lim_{n \rightarrow \infty} \phat_n(p_n) =  
   \phi(F(p))    
\]
when $p_n \longrightarrow p$ with $ p_n \in \e_n(S^n)$.  
We are assuming $S$ is complete and separable, therefore so is 
$\P(S)$ \cite{Dudley}.

We may now apply Lemma \ref{lemma} with 
\[
\begin{array}{ccc}
      X = \P(S),    &   D_n = \e_n(S^n), &  f_n = \phat_n,   \\
\end{array}
\]
and conclude that 
\begin{equation}
\label{eq}
  \lim_{n \rightarrow \infty} \int_{\P(S)} \phat_n(\zeta) \mu_n(d\zeta) =
          \phi(F(p))  
\end{equation}
for any sequence $\{ \mu_n  \}$ that converges to 
$\delta(p)$ in $\P(\P(S))$.  

By equations (\ref{eq}) and (\ref{phat}), 
\begin{eqnarray*}
 \phi(F(p)) &=&
 \lim_{n \rightarrow \infty} \int_{\P(S)} \phat_n(\zeta) \mu_n(d\zeta)  \\ 
 &=&
 \lim_{n\rightarrow \infty}\int_{\P(S)}\int_{\P(T)} \phi(\eta) H_n(\zeta,d\eta)
		    \mu_n(d\zeta)                      \\ 
 &=& \lim_{n \rightarrow \infty} \int_{\P(T)} \phi(\eta) 
	    \int_{\P(S)} H_n(\zeta, d\eta) \mu_n(d\zeta)    ,  \\ 
\end{eqnarray*}
for all $\phi \in C_b(\P(T))$.
This implies that 
$
    \int_{\P(S)} H_n(\zeta, d\eta) \mu_n(d\zeta) \longrightarrow 
    \delta(F(p))     
$
in $\P(\P(S))$, completing the proof.
\qquad $\blacksquare$

\section{Acknowledgements}

I would like to thank A. J. Chorin for introducing me to the
topic of propagation of chaos, and for referring me to a book \cite{Kac} of M. Kac 
and to research of D. Talay and M. Bossy \cite{BossyTalay}.  I am very
grateful to L. Le Cam for countless helpful conversations and for lending me {\it
Convergence of Probability Measures} by P. Billingsley
\cite{Billingsley}.

\flushright

MATHEMATICS DEPARTMENT 

COMPUTING SCIENCES DIRECTORATE

LAWRENCE BERKELEY NATIONAL LAB

1 CYCLOTRON ROAD, BERKELEY, CA  94720

\end{document}